\documentclass{birkjour}
\usepackage{amsthm}
\usepackage{amsfonts}
\usepackage{mathrsfs}
\usepackage{amscd,xy,amssymb,amsmath,graphicx,verbatim}
\usepackage[bookmarksnumbered, colorlinks, plainpages]{hyperref}

\newtheorem{theorem}{Theorem}[section]
\newtheorem{lemma}[theorem]{Lemma}
\newtheorem{corollary}[theorem]{Corollary}
\newtheorem{proposition}[theorem]{Proposition}

\newtheorem{definition}[theorem]{Definition}
\newtheorem{example}[theorem]{Example}

\theoremstyle{remark}
\newtheorem{remark}[theorem]{Remark}

\numberwithin{equation}{section}

\begin{document}

%
%
%
%
%
%
%
%
%

\title{Topological Entropy Conjecture}

\author[Lvlin Luo]{Lvlin Luo}
\address{Mathematics School, Jilin University, 130012, Changchun, P. R. China}
\address{School of Mathematical Sciences, Fudan University, 200433, Shanghai, P. R. China}
\address{School of Mathematics and Statistics, Xidian University, 710071, Xi'an, P. R. China}
\email{luoll12@mails.jlu.edu.cn}
\email{luolvlin@fudan.edu.cn}

\subjclass{Primary 37B40, 54H20; Secondary 28D20, 57T99.}

\keywords{algebra equation; \v{C}ech homology; \v{C}ech homology germ; eigenvalue; topological fiber entropy.}

\begin{abstract}
In 1974, M. Shub stated Topological Entropy Conjecture,
that is, the inequality $\log\rho\leq ent(f)$ is valid or not,
where $f$ is a continuous self-map on a compact manifold $M$,
$ent(f)$ is the topological entropy of $f$ and
$\rho$ is the maximum absolute eigenvalue of $f_*$
which is the linear transformation induced by $f$ on the homology group $H_{*}(M;\mathbb{Z})=\bigoplus\limits_{i=0}^n{H_{i}(M;\mathbb{Z})}$.
In 1986, A. B. Katok gave a counterexample such that the inequality $\log\rho\leq ent(f)$ is invalid.
In this paper, we define $f$-\v{C}ech homology group $\check{H}_{i}(X,f;\mathbb{Z})$
and topological fiber entropy $ent(f_L)$ on compact Hausdorff space $X$
for which there is $n=n(J)$ such that $\check{H}_*(X;\mathbb{Z})$ exists,
where $f\in C^0(X)$ and $J$ is the set of all covers.
Then we prove that $\log\rho\leq ent(f_L)$ is valid.
\end{abstract}

\maketitle

\section{Introduction}

In 1974, M. Shub stated a conjecture,
named Topological Entropy Conjecture \cite{MShub1974}.
That is,
let $f\in C^0(M)$,
then $f$ induces a homomorphism
$f_*:$ $H_{*}(M;\mathbb{Z})\to H_{*}(M;\mathbb{Z})$,
where $\mathbb{Z}$ is the set of integer,
$M$ is a $n$-dimension compact manifold,
$C^0(M)$ is the set of all continuous self-mapping on $M$,
$H_{*}(M;\mathbb{Z})=\bigoplus\limits_{i=0}^n{H_{i}(M;\mathbb{Z})}$
and $H_{i}(M;\mathbb{Z})$ is the $i$th integer coefficients homology group of $M$.
Hence,
on the homology group $H_{*}(M;\mathbb{Z})$,
$f_{*}$ is a linear transformation which is a $(n+1)\times(n+1)$ matrix concerning integer entries.
Let $\rho$ be the maximum absolute eigenvalue of $f_{*}$,
then the topological entropy conjecture is that the inequality $\log\rho\leq ent(f)$ is valid or not,
where $ent(f)$ is the topological entropy of $f$.

The inequality is so simple connected in the first place with the work of Smale, Shub, and Sullivan,
that one attempts to prove it have been very fruitful.
Newhouse and Yomdin proved this for $C^\infty$ diffeomorphisms.
But unlucky,
in 1986,
A. B. Katok gave a counterexample explaining the inequality invalid for continuous maps,
that is,
on $2$-dimension sphere $S^{2}$ there is $f\in C^{0}(S^2)$ such that
$0=ent(f)<\log \rho$ \cite{ABKatok1986}.

After studying the counterexample,
the definition of topological entropy and $\check{C}ech$ cohomology,
through a little revision in definition,
we proof the conjecture is valid again with our definition
on compact Hausdorff space $X$
for which there is $n=n(J)$ such that $\check{H}_*(X;\mathbb{Z})$ exists,
where $f\in C^0(X)$ and $J$ is the set of all covers.
To establish the inequality is my interesting,
inevitable there is a lot of words about homology theory in this paper.

To be short of reading the origin papers,
there is a story between Klein and Poincar\`{e} about Fuchs function in $2$-dimension,
of course, in the last Poincar\`{e} was remedial this by named Klein group from his own achievement in $3$-dimension.
The same thing,
some times I maybe forgot or ignore the references for never reading the origin papers,
if someone who find something, please forgive me and chase me, I will remedy that in the first time.

Generally, in this paper,
$X$ denotes a compact Hausdorff space,
$C^0(X)$ denotes the set of all continuous self-mapping on $X$,
$id$ or $I$ denotes the identical mapping on $X$ and
$\mathbb{Z}$ is the integer group.
The set of covers always is the set of open covers.
For brevity,
$n$ and $\|\cdot\|$ denotes many kinds of meaning,
one can regard this following the context.

Let $\alpha,\beta$ be open covers of $X$,
if for any $B\in\beta$,
exist $A\in\alpha$ such that $B\subseteq A$,
then define $\alpha<\beta$,
named that $\beta$ is larger than $\alpha$ or $\beta$ is a refinement of $\alpha$.

Let $\alpha^c=\{A; A^c\in\alpha,\text{where } A^c\cup A=X \text{and } A^c\cap A=\emptyset\}$,
for $A\in\alpha$,
$\overline{A}$ denotes the closure of $A$,
$\|A\|$ denotes the numbers of elements of $A$.

Denotes:
$$\left\{
   \begin{array}{lr}
    a_{0}\cdots a_{i-1}\hat a_{i},a_{i+1}\cdots a_{p}=a_{0}\cdots a_{i-1}a_{i+1},\cdots a_{p};\\
    a_{0}\cdots a_{i-1}\underline b^{(i)} a_{i}\cdots a_{p}=a_{0}\cdots a_{i-1}ba_{i}\cdots a_{p};\\
    a_{0}\cdots a_{i-1}\underline b_{(k)}^{(i)} a_{i}\cdots a_{p}=\sum\limits_{m\in (k)}a_{0}\cdots a_{i-1}b_{m}a_{i}\cdots a_{p};\\
    a_{0}\cdots a_{i-1}\underline b_{\emptyset}^{(i)}\cdots a_{p}=\sum\limits_{m\in \emptyset}a_{0}\cdots a_{i-1}b_{m}a_{i}\cdots a_{p}=a_{0}\cdots a_{i-1}a_{i}\cdots a_{p};\\
    (k)=\{k_1,k_2,k_3,\cdots,k_n;n=\|\{a_{0},\cdots,a_{i-1},b_{m},a_{i},\cdots,a_{p}\}\|\geq1,m\in\mathbb{Z}\};\\
    {(a_{0}\cdots\hat{b_{(k)}}\cdots a_{p})}^d=b_{k_1}\cdots b_{k_i}\cdots b_{k_n},k_i\in(k).
    \end{array}
 \right.$$
Before starting the main body, first please me express gratitude to Professor Bingzhe Hou for helpful suggestions on Definition~\ref{3}.
I respect my PhD supervisor Youqing Ji,
for his selfless and the spirit of education.
Finally thanks
Professor Yang Cao, classmates, friends, and Jilin University,
in the one or two years after my father died, that was my cloudy day,
for their help in my daily lives I life better.

\section{Algebra Equation for Boundary Operator}

\begin{definition}[\cite{JamesRMunkres2006}P541]\label{1}
Let $X$ be a Hausdorff space and
$\Psi$ be a cover of $X$ and $U_0,U_1,U_2,\cdots,U_p\in\Psi$,
if $U_0\bigcap U_{1}\bigcap\cdots\bigcap U_{p}\neq\emptyset$,
then define a $p$-simplex $\sigma_{p}$ and $p$th chain group $C_{p}$,
so get the $p$th homology group $H_{p}(\Psi;\mathbb{Z})$ and cohomology group $H^{p}(\Psi;\mathbb{Z})$,

where
$$\cdots\xrightarrow{}C_{p+1}(\Psi;\mathbb{Z})\xrightarrow{\partial_{p+1}}C_{p}(\Psi;\mathbb{Z})\xrightarrow{\partial_{p}}C_{p-1}(\Psi;\mathbb{Z})\xrightarrow{}\cdots$$

$$\partial_{p} (U_{0}\cap\cdots\cap U_{p})=\sum\limits_{i}^{p}(-)^{i}(U_{0}\cap\cdots\cap \hat{U_{i}}\cdots\cap U_{p}),$$

$$\partial_{p-1}\circ\partial_{p}=0,$$

$$\left.
   \begin{array}{lr}
    B_{p}(\Psi;\mathbb{Z})=im\partial_{p+1},\\
    Z_{p}(\Psi;\mathbb{Z})=\ker\partial_{p},\\
    {H}_{p}(\Psi;\mathbb{Z})=Z_{p}/B_{p}.\\
    \end{array}
 \right.$$

Let $$C^{p}(\Psi;\mathbb{Z})=Hom(C_{p}(\Psi;\mathbb{Z}),\mathbb{Z}),$$

so $\partial_{p}$ induce a homomorphism $C^{p-1}(\Psi;\mathbb{Z})\xrightarrow{\delta^{p}} C^{p}(\Psi;\mathbb{Z})$,

there is
$$\cdots\xleftarrow{}C^{p+1}(\Psi;\mathbb{Z})\xleftarrow{\delta^{p+1}}C^{p}(\Psi;\mathbb{Z})\xleftarrow{\delta^{p}}C^{p-1}(\Psi;\mathbb{Z})\xleftarrow{}\cdots,$$

$$   \delta^{p+1}\circ\delta^{p}=0,$$

$$\left.
   \begin{array}{lr}
    B^{p}(\Psi;\mathbb{Z})=im\delta^{p},\\
    Z^{p}(\Psi;\mathbb{Z})=ker\delta^{p+1},\\
    {H}^{p}(\Psi;\mathbb{Z})=Z^{p}/B^{p}.\\
    \end{array}
 \right.$$

\end{definition}

\begin{lemma}\label{2}
Let $X$ be a Hausdorff space and
$\Psi$ be a cover of $X$,
we get $C^{p}(\Psi;\mathbb{Z})\cong C_{p}(\Psi;\mathbb{Z})$.
Moreover, let $U^{i}=U_{i}^c$,
if $c_p=U_{0}\cap\cdots\cap U_{p}\in C_{p}(\Psi;\mathbb{Z})$,
then $c^p=U^{0}\cup\cdots\cup U^{p}\neq X$ is a isomorphic representation of $p$-simplex of $C^{p}(\Psi;\mathbb{Z})$.
\end{lemma}
\begin{proof}
Because $\mathbb{Z}$ can be regard as a finite generated free ring \cite{FrankWAndersonKentRFuller1992},
$C_{p}(\Psi;\mathbb{Z})$ can be treated as a finite dimension $\mathbb{Z}$-module space \cite{BLvanderWaerden1991},
and $C^{p}(\Psi;\mathbb{Z})$ can be treated as the dual $\mathbb{Z}$-module space of $C_{p}(\Psi;\mathbb{Z})$.
With the property of finite dimension $\mathbb{Z}$-module space,
we get $C^{p}(\Psi;\mathbb{Z})\cong C_{p}(\Psi;\mathbb{Z})$.

$$c_p=U_{0}\cap\cdots\cap U_{p}\neq\emptyset\Longleftrightarrow c^p=U^{0}\cup\cdots\cup U^{p}\neq X,$$

so

$$c_p\in C_{p}(\Psi;\mathbb{Z})\Longleftrightarrow c^p\in C^{p}(\Psi;\mathbb{Z}),$$
that is,
$U^{0}\cup\cdots\cup U^{p}\neq X$ is a isomorphic representation of $p$-simplex of $C^{p}(\Psi;\mathbb{Z})$.
\end{proof}

\begin{definition}\label{3}
Let $X$ be a Hausdorff space
and $\Psi$ be a cover of $X$,
then define
$n_\Psi=\max\{n;n\in S\}$
and named that $n_\Psi$ is the $\partial$ operator dimension of $X$ about $\Psi$
where
\begin{small}
$$S=\{\partial(U_0\cdots\cap U_i\cdots U_n)\neq\partial(U_0\cdots\cap U_i\cdots U_n\cap U_{n+1});U_0,\cdots,U_{n+1}\in\Psi\}.$$
\end{small}
Let $J$ be the ordered set induced by the refinement cover of $X$,
obviously, if $\alpha,\beta\in J,\alpha<\beta$,
then $n_\alpha\leq n_\beta$.
So $n_\Psi=n(\Psi)$ is a function defined on $J$.
If $\lim\limits_{\overrightarrow{\Psi\in J}} n_\Psi=n_J=n(J)$ exist,
then define $n_J$ is the $\partial$ operator dimension of $X$.
\end{definition}

In order to deal with the problem easily,
always let $\Psi\in J$ be good enough and enough refinement,
i.e.,
satisfy all the necessary requirements of the problem.

\begin{definition}\label{4}
Let $X$ be a Hausdorff space
and $\Psi$ be a cover of $X$ and $p\leq n=n(\Psi)$.
If for any $\sigma^p\in C^{p}(\Psi;\mathbb{Z})$,
there is $\sigma^n\in C^{n}(\Psi;\mathbb{Z})$ such that $\sigma^p=U^{0}\cup\cdots\cup U^{p}$ is the $p$th surface of $\sigma^n$
and there is
$${(U^{0}\cup\cdots\cup \hat{U_{(k)}}\cdots\cup U^{p})}^d=U^{k_0}\cup\cdots\cup U^{k_{n-p+1}}.$$
Then define $X$ is a Poincar\`{e} space.
\end{definition}

\begin{lemma}\label{5}
Let $X$ be a Poincar\`{e} space and $\Psi$ be a cover of $X$,
there is
${H}^{p}(\Psi;\mathbb{Z})\cong {H}_{n-p}(\Psi;\mathbb{Z})$.
\end{lemma}
\begin{proof}

By $lemma~\ref{2}$ get two chains:

\hspace*{-1cm}(1)\begin{small}
$\left\{
   \begin{array}{lr}
    \cdots\xrightarrow{}C_{p+1}(\Psi;\mathbb{Z})\xrightarrow{\partial_{p+1}}C_{p}(\Psi;\mathbb{Z})\xrightarrow{\partial_{p}}C_{p-1}(\Psi;\mathbb{Z})\xrightarrow{}\cdots\\
    \cdots\xleftarrow{}C^{p+1}(\Psi;\mathbb{Z})\xleftarrow{\delta^{p+1}}C^{p}(\Psi;\mathbb{Z})\xleftarrow{\delta^{p}}C^{p-1}(\Psi;\mathbb{Z})\xleftarrow{}\cdots\\
    \end{array}
 \right.$
 \end{small}

for a fixed $p$-simplex of $C_{p}(\Psi;\mathbb{Z})$,considered the algebra equation:

\hspace*{-1cm}(a)\begin{small}
$\left\{
   \begin{array}{lr}
    <\partial c_p,c^{p-1}>=<c_p,\delta c^{p-1}>,\\
    \partial_{p} (U_{0}\cdots\cap U_i\cdots U_{p})=\sum\limits_{i=0}^{p}{(-1)}^{i}(U_{0}\cdots\cap \hat{U_{i}}\cdots U_{p}),\\
    \partial\emptyset=\delta\emptyset=0,\\
    <a,\emptyset>=<\emptyset,b>=0.
    \end{array}
 \right.$
 \end{small}

\hspace*{-1cm}(b)\begin{small}
$\left\{
   \begin{array}{lr}
    <\sum\limits_{i=0}^{p}{(-1)}^{i}(U_{0}\cdots\cap \hat{U_{i}}\cdots U_{p}),V_{0}\cdots\cap \hat{V_{i}}\cdots V_{p}>=<c_p,\delta c^{p-1}>,\\
    \sum\limits_{i=0}^{p}(-)^{i}(V^{0}\cdots\cup \underline{V_{(k)}^{(i)}}\cdots V^{p})=\delta^{p}(V^{0}\cdots\cup\hat{V_{(k)}}\cdots V^{p}),\\
    \end{array}
 \right.$
 \end{small}

if $(k)=\emptyset$,
then define $$(i)=\emptyset,{(-1)}^{\emptyset}=0,$$
so $$\delta^{p}(U^{0}\cdots\cup\hat{U_{(k)}}\cdots U^{p})=0.$$
From (a)(b), we get

\hspace*{-1cm}(2)\begin{small}
$\left\{
   \begin{array}{lr}
    \partial_{p} (U_{0}\cdots\cap U_i\cdots U_{p})=\sum\limits_{i=0}^{p}{(-1)}^{i}(U_{0}\cdots\cap \hat{U_{i}}\cdots U_{p}),\\
    \delta^{p}(U^{0}\cdots\cup \hat{U_{(k)}}\cdots U^{p})=\sum\limits_{i=0}^{p}{(-1)}^{i}(U^{0}\cdots\cup \underline{U_{(k)}^{(i)}}\cdots U^{p}).\\
    \end{array}
 \right.$
 \end{small}

that is :

\hspace*{-1cm}(3)\begin{small}
$\left\{
   \begin{array}{lr}
    \partial_{p}(U_{0}\cdots\cap U_i\cdots U_{p})-\sum\limits_{i=0}^{p}{(-1)}^{i}(U_{0}\cdots\cap \hat{U_{i}}\cdots U_{p})=0,\\
    \delta^{p}(U^{0}\cdots\cup \hat{U_{(k)}}\cdots U^{p})-\sum\limits_{i=0}^{p}{(-1)}^{i}(U^{0}\cdots\cup \underline{U_{(k)}^{(i)}}\cdots U^{p})=0,\\
     \delta^{n-p+1}({(U^{0}\cdots\cup \hat{U_{(k)}}\cdots U^{p})}^d)=\delta^{n-p+1}{(U^{k_1}\cdots\cup {U^{k_m}}\cdots U^{k_{n-p+1}})},\\
     \delta^{n-p+1}{(U^{k_1}\cdots\cup {U^{k_m}}\cdots U^{k_{n-p+1}})}=\sum\limits_{i=0}^{p}{(-1)}^{i}{(U^{k_1}\cdots\cup \underline{U_{(0,\cdots, p)}^{(i)}}\cdots U^{k_{n-p+1}})}.\\
    \end{array}
 \right.$
 \end{small}

Let $$c_p=\sum z_m{(U_{0}\cdots\cap U_i\cdots U_{p})}_m ,$$

so $$c^{n-p}=\sum z_m{({(U^{0}\cdots\cup \hat{U_{(k)}}\cdots U^{p})}^d)}_m,$$
where $z_m\in\mathbb{Z}.$

\hspace*{-1cm}(c)\begin{small}
$\left\{
   \begin{array}{lr}
U_{0}\cdots\cap U_i\cdots U_{p}\longleftrightarrow U^{0}\cdots\cup \hat{U_{(k)}}\cdots U^{p}\longleftrightarrow {(U^{0}\cdots\cup \hat{U_{(k)}}\cdots U^{p})}^d,\\
c_p\in ker\partial_{p}\Longleftrightarrow
c^{n-p}\in ker\delta^{n-p+1},\\
c_p \in im\partial_{p+1}\Longleftrightarrow
c^{n-p}\in im\delta^{n-p}.
   \end{array}
 \right.$
 \end{small}

Let:
(4)\begin{small}
$\left\{
   \begin{array}{lr}
   \partial_{\frac{ker}{im}}(C_p)={H}_{p}(\Psi;\mathbb{Z})=Z_{p}/B_{p}=ker\partial_{p}/im\partial_{p+1},\\
   \partial^*_{\frac{ker}{im}}(C^p)={H}^{p}(\Psi;\mathbb{Z})=Z^{p}/B^{p}=ker\delta^{p+1}/im\delta^{p}.\\
   \end{array}
 \right.$
 \end{small}

So $\partial_{p}$ and $\delta^{n-p+1}$ is dual solution in algebra equation (3),
similarly $\partial_{\frac{ker}{im}}$ and $\partial^*_{\frac{ker}{im}}$ is dual value in (4),
the all process of dual mapping is linear reversible,
i.e., the same style as isomorphism.
Therefore,
the $p$th value of $\partial_{\frac{ker}{im}}$ on the $C_p$ chain group is isomorphic to the $(n-p)$th value of $\partial^*_{\frac{ker}{im}}$ on the $C^{n-p}$ chain group,
that is $\partial_{\frac{ker}{im}}(C_p)\cong\partial^*_{\frac{ker}{im}}(C^{n-p})$,
for this reason, there is
$${H}^{p}(\Psi;\mathbb{Z})\cong {H}_{n-p}(\Psi;\mathbb{Z}).$$
\end{proof}

Likely linear equation of Euclidean space $\mathbb{R}^3$,
let $S_i:A_ix+B_iy+C_iz=0$ be a class of lines,
or in another word that is a class of planes $S^*_i:A_ix+B_iy+C_iz=0$,
where $0\leq i<+\infty$.

Line and plane is a pair of dual,
and for a fixed space,
the intrinsic relationships between line or between plane are never changed.
That is, let $f,g$ be good mappings,
such as there are linear,
if $f_i=f(S_i,S_{i-1}),f^*_i=f^*(S^*_i,S^*_{i+1})$
and $g_i=g(f_i),g^*_i=g(f^*_i)$,
then $g_i,g^*_i$ is a pair of dual,
such that there is a natural relationship between $g_i$ and $g^*_{n-i}$.
For example,
that natural relationship maybe is $g_i=g^*_{n-i}$,
or $g_ig^*_{n-i}=1$,
or $g_i+g^*_{n-i}=0$,
or $g_iA_k+g_{n-i}B_k+C_k=0$ and $g_{n-i}^{*}A_k+g_{i}^{*}B_k+C_k=0$,
and so on,
the dual outcomes and the representations of the natural relation between $g_i$ and $g_{n-i}^{*}$ only dependent on the good mappings $f,g$.

\section{Germ and Dual of \v{C}ech homology}

\begin{definition}[\cite{JamesRMunkres2006}P542]\label{6}
Let $X$ be a Hausdorff space
and $J$ be the ordered set induced by the set of all covers of $X$,
$U_0,U_1,U_2,\cdots,U_p\in\Psi,\Psi\in J$,
if $U_0\bigcap U_{1}\cap\cdots\cap U_{p}\neq\emptyset$,
then define a $p$-simplex $\sigma_{p}$ and we get the $p$th chain group $C_{p}$,
moreover we get the $p$th homology group $H_{p}(\Psi;\mathbb{Z})$ and the $p$th cohomology group $H^{p}(\Psi;\mathbb{Z})$.
If $\Omega,\Psi\in J$ and $\Omega<\Psi$,
then we get the homomorphisms
$$f_{\Psi\Omega}:H_{p}(\Psi;\mathbb{Z})\to H_{p}(\Omega;\mathbb{Z}),$$
$$f_{\Omega\Psi}:H^{p}(\Omega;\mathbb{Z})\to H^{p}(\Psi;\mathbb{Z}).$$
then define the
$p$th \v{C}ech cohomology group
$$\check{H}^{p}(X;\mathbb{Z})=\lim\limits_{\overrightarrow{\Omega\in J}}H^{p}(\Omega;\mathbb{Z}).$$
\end{definition}

\begin{definition}\label{9}
Following Definition~$\ref{6}$,
define the $p$th \v{C}ech homology group
$$\check{H}_{p}(X;\mathbb{Z})=\lim\limits_{\overleftarrow{\Omega\in J}}H_{p}(\Omega;\mathbb{Z}).$$
\end{definition}

\begin{definition}\label{7}
Let $X$ be a Poincar\`{e} space
and $J$ be the ordered set induced by the set of all covers of $X$,
$n=n(J)$.
For $\Omega,\Psi\in J$, let $\Theta=\Psi\vee\Omega=\{\alpha\cap\beta;\alpha\in\Psi,\beta\in\Omega\}$,
then with Definition~$\ref{6}$,
we get homomorphisms $f_{\Theta\Omega}:H_{p}(\Theta;\mathbb{Z})\to H_{p}(\Omega;\mathbb{Z})$
and $f_{\Theta\Psi}:H_{p}(\Theta;\mathbb{Z})\to H_{p}(\Psi;\mathbb{Z})$.
Following this, we can define a \v{C}ech homology germ $H_p(J;\mathbb{Z})$,
and similarly define a \v{C}ech cohomology germ $H^p(J;\mathbb{Z})$.
If exist $\Gamma\in J$, for any $\Psi\in J$ when $\Gamma<\Psi$,
there is $H_{p}(\Psi;\mathbb{Z})\cong H^{n-p}(\Psi;\mathbb{Z})$,
then define
$$H^p(J;\mathbb{Z})\cong H_{n-p}(J;\mathbb{Z}).$$
\end{definition}

\begin{lemma}\label{10}
Let $X$ be a Hausdorff space
$X$
for which there is $n=n(J)$ such that $\check{H}_p(X;\mathbb{Z})$ exist,
where $p\leq n$ and $J$ is a ordered set induced by the set of all covers of $X$,
there are
$\check{H}_{p}(X;\mathbb{Z})\sim H_{p}(J;\mathbb{Z})$ and
$\check{H}^{p}(X;\mathbb{Z})\sim H^{p}(J;\mathbb{Z})$,
where $\sim$ means the different expressions of the same thing.
\end{lemma}
\begin{proof}
By Lemma~$\ref{5}$ and Definition~$\ref{6},\ref{7},\ref{9}$.
\end{proof}

\begin{definition}\label{11}
Let $X$ be a Poincar\`{e} space
and $J$ be a ordered set induced by the set of all covers of $X$,
$n=n(J)$.
If $H_{p}(J;\mathbb{Z})\cong H^{n-p}(J;\mathbb{Z})$,
then define $\check{H}_{p}(X;\mathbb{Z})\cong\check H^{n-p}(X;\mathbb{Z})$.
\end{definition}

\section{f-\v{C}ech homology}

\begin{definition}\label{13}
Let $X$ be a Hausdorff space,
$U_i,V,W\subseteq X$ and $f\in C^0(X)$,
where $0\leq i\leq k,k\in \mathbb{Z}$.
define:

\hspace*{-1cm}\begin{small}
$\left\{
   \begin{array}{l}
    L_{f}(U)=\{\cdots,{f}^{-n}(U),\cdots,{f}^{-1}(U),{f}^{0}(U),\cdots,{f}^{n}(U),\cdots\};\\
    f\circ L_f=L_f\circ f;\\
    L_{f}(U)\cap L_{f}(V)=L_{f}(W),\text{where } W=U\cap V;\\
    L_{f}(U_0)\cdots \cap L_{f}(U_i)\cdots L_{f}(U_k)=L_{f}(U_0)\cap(L_{f}(U_1)\cdots\cap L_{f}(U_i)\cap L_{f}(U_k));\\
    L_{g+h}(U)=\{\cdots,{g}^{-n}(U)\bigcup{h}^{-n}(U),\cdots,{g}^{0}(U)\bigcup {h}^{0}(U),\cdots,{g}^{n}(U)\bigcup{h}^{n}(U),\cdots\};\\
    L_{f}(\emptyset)=\emptyset;\\
    L_{g\oplus h}(U)=L_{g+h}(U),when\,{g}^{-1}(U)\cap{h}^{-1}(U)=\emptyset.
    \end{array}
 \right.$
 \end{small}

then $L_{f}(U)$ is named the $f$-fiber of $U$,
and define
$$X^{f}=\{L_{f}(U);U\subset X\}.$$
\end{definition}

If $X$ is a compact space,
then
$X^{+\infty}=\prod\limits_{i=-\infty}^{-1} X\times\underline{X}\times\prod\limits_{i=1}^{+\infty} X$
is compact too by Tychonoff theorem.

\begin{definition}\label{14}
Let $X$ be a Hausdorff space and let $J$ be the ordered set induced by the set of all covers of $X$,
$f\in C^0(X),\Psi\in J$,
$U_0,\cdots,U_p\in\Psi$.
If $\sigma^f_{p}=L_{f}(U_{0})\cap\cdots\cap L_{f}(U_{p})\neq\emptyset$,
then define a $f$-\v{C}ech $p$-simplex $\sigma^f_{p}$
and the $f$-\v{C}ech $p$-chain group ${C}_{p}(\Psi,f;\mathbb{Z})$,
so get the $f$-\v{C}ech $p$th homology group ${H}_{p}(\Psi,f;\mathbb{Z})$,
where

\hspace*{-1cm}\begin{small}
$\left.
   \begin{array}{c}
\cdots\xrightarrow{}{C}_{p+1}(\Psi,f;\mathbb{Z})\xrightarrow{\partial^f_{p+1}}{C}_{p}(\Psi,f;\mathbb{Z})\xrightarrow{\partial^f_{p}}{C}_{p-1}(\Psi,f;\mathbb{Z})\xrightarrow{}\cdots,\\
\partial^f_{p} (L_f(U_{0})\cdots\cap L_f(U_{i})\cdots L_f(U_{p}))=\sum\limits_{i=0}^{p}{(-1)}^{i}(L_f(U_{0})\cdots\cap \hat{L}_f(U_{i})\cdots L_f(U_{p})),
    \end{array}
 \right.$
 \end{small}

It is easy to get $$\partial^f\circ\partial^f=0,$$
in fact, there is

\begin{small}
$\left.
   \begin{array}{l}
\partial^f\circ\partial^f (L_f(U_{0})\cdots\cap L_f(U_{i})\cdots L_f(U_{p}))\\
=\sum\limits_{i}^{p}{(-1)}^{i}\partial^f(L_f(U_{0})\cdots\cap\hat{L}_f(U_{i})\cdots L_f(U_{p}))\\
=\sum\limits_{i}^{p}\sum\limits_{j<i}{(-1)}^{i+j}(L_f(U_{0})\cdots\cap \hat{L}_f(U_{j})\cdots\cap \hat{L}_f(U_{i})\cdots L_f(U_{p}))\\
+\sum\limits_{i}^{p}\sum\limits_{j>i}{(-1)}^{i+j-1}(L_f(U_{0})\cdots\cap \hat{L}_f(U_{i})\cdots\cap \hat{L}_f(U_{j})\cdots\cap L_f(U_{p}))\\
=0,
    \end{array}
 \right.$
 \end{small}

and
$$\left.
   \begin{array}{l}
    {B}_{p}(\Psi,f;\mathbb{Z})=im\partial^f_{p+1},\\
    {Z}_{p}(\Psi,f;\mathbb{Z})=ker\partial^f_{p},\\
    {H}_{p}(\Psi,f;\mathbb{Z})=Z_{p}(\Psi,f;G)/B_{p}(\Psi,f;\mathbb{Z}),\\
    \end{array}
 \right.$$

\end{definition}

\begin{lemma}\label{15}
A \v{C}ech $p$-chain $c_{p}$ can induce a $f$-\v{C}ech $p$-chain $c^f_{p}$,
that is, $U_0\cap U_1\cap\cdots\cap U_p\neq\emptyset$ if and only if
$L_f(U_0)\cap L_f(U_1)\cap\cdots\cap L_f(U_p)\neq\emptyset$.
So the \v{C}ech $p$-chain group is isomorphic to the $f$-\v{C}ech $p$-chain group.
\end{lemma}
\begin{proof}
Using Lemma~$\ref{10}$,definition~$\ref{13}$.
\end{proof}

\begin{definition}\label{3ftongdiaoweishudingyi}
Let $X$ be a Hausdorff space and $f\in C^0(X)$ and
$\Psi$ be a cover of $X$,
then define
$n_{\Psi,f}=\max\{n;n\in S\}$
and named that $n_{\Psi^f}$ is the $\partial^f$ operator dimension of $X$ about $\Psi$,
where

\begin{small}
$$S=\{\partial^f(L_f(U_0)\cdots \cap L_f(U_n))\neq\partial(L_f(U_0)\cdots \cap L_f(U_n)\cap L_f(U_{n+1}));U_0,\cdots,U_{n+1}\in\Psi\}.$$
\end{small}

Let $J$ be the ordered set induced by the refinement cover of $X$,
obviously, if $\alpha,\beta\in J,\alpha<\beta$,
then $n_{\alpha,f}\leq n_{\beta,f}$.
So $n_{\Psi,f}=n(\Psi,f)$ is a function defined on $J$ for fixed $f$.
If $\lim\limits_{\overrightarrow{\Psi\in J}} n_{\Psi,f}=n_{J,f}=n(J,f)$ exist,
then define $n_{J,f}$ is the $\partial^f$ operator dimension of $X$.
\end{definition}

\begin{definition}\label{17}
Following definition$~\ref{3ftongdiaoweishudingyi}$,
similarly with definition$~\ref{7},~\ref{9}$,
Let $X$ be a Hausdorff space,
$f\in C^0(X)$ and
$J$ be the ordered set induced by the refinement cover of $X$,
$n=n(J,f)$.
For $\Omega,\Psi\in J$, let $\Theta=\Psi\vee\Omega=\{\alpha\cap\beta;\alpha\in\Psi,\beta\in\Omega\}$,
then similarly with definition~$\ref{6}$,
we get homomorphisms $f_{\Theta\Omega}:H_{p}(\Theta,f;\mathbb{Z})\to H_{p}(\Omega,f;\mathbb{Z})$
and $f_{\Theta\Psi}:H_{p}(\Theta,f;\mathbb{Z})\to H_{p}(\Psi,f;\mathbb{Z})$.
Following this, similarly with definition$~\ref{7}$,
we can define a $p$th $f$-\v{C}ech homology germ $H_p(J,f;\mathbb{Z})$.
And similarly with definition$~\ref{9}$,
we can define the $p$th $f$-\v{C}ech homology group
$$\check{H}_{p}(X,f;\mathbb{Z})=\lim\limits_{\overleftarrow{\Omega\in J}}H_{p}(\Omega,f;\mathbb{Z}).$$
\end{definition}

\begin{lemma}\label{ffibertongdiaoqunyukongjiantongdiaoquntonggou1}
Let $X$ be a Hausdorff space and $f\in C^0(X)$
for which there is $n=n(J,f)$ such that $\check{H}_p(X,f;\mathbb{Z})$ exist,
where $p\leq n$ and $J$ is a ordered set induced by the set of all covers of $X$.
Then there is $n=n(J)$ such that $\check{H}_p(X;\mathbb{Z})$ exist for $p\leq n$,
moreover
$$\left.
   \begin{array}{c}
    n=n(J)=n(J,f),\\
    im\partial_{p+1}={B}_{p}(\Psi;\mathbb{Z})={B}_{p}(\Psi,f;\mathbb{Z})=im\partial^f_{p+1},\\
    ker\partial_{p}={Z}_{p}(\Psi;\mathbb{Z})={Z}_{p}(\Psi,f;\mathbb{Z})=ker\partial^f_{p},\\
    Z_{p}(\Psi;G)/B_{p}(\Psi;\mathbb{Z})={H}_{p}(\Psi;\mathbb{Z})={H}_{p}(\Psi,f;\mathbb{Z})=Z_{p}(\Psi,f;G)/B_{p}(\Psi,f;\mathbb{Z}),\\
    \end{array}
 \right.$$
\end{lemma}

\begin{lemma}\label{18}
Let $X$ be a Hausdorff space and $f\in C^0(X)$
for which there is $n=n(J,f)$ such that $H_p(X,f;\mathbb{Z})$ exist,
where $p\leq n$ and $J$ is a ordered set induced by the set of all covers of $X$,
then there is
$$H_{p}(X,J;\mathbb{Z})\sim \check{H}_{p}(X,f;\mathbb{Z}),$$
where $\sim$ means the different expressions of the same thing.
\end{lemma}
\begin{proof}
Using lemma~$\ref{10}$,
definition~$\ref{17}$,
lemma~$\ref{ffibertongdiaoqunyukongjiantongdiaoquntonggou1}$.
\end{proof}

Also,
we can define $f$-\v{C}ech cohomology germ $H^{p}(J,f;\mathbb{Z})$,
$f$-\v{C}ech cohomology group $\check{H}^{p}(X,f;\mathbb{Z})$
and $f$-Poincar\`{e} space.
From the definition, there is obviously
$${C}_{p}(X;\mathbb{Z})={C}_{p}(X,id;\mathbb{Z})$$.
For convenience define:

\begin{small}$$\left.
   \begin{array}{l}
    \check{H}_{*}(X;\mathbb{Z})=\bigoplus\limits_{i=0}^n{\check H_{i}(X;\mathbb{Z})};\\
    C_{*}(X;\mathbb{Z})=\bigoplus\limits_{i=0}^n{C_{i}(X;\mathbb{Z})};\\
    B_{*}(X;\mathbb{Z})=\bigoplus\limits_{i=0}^n{B_{i}(X;\mathbb{Z})};\\
    \check{H}_{*}(X,f;\mathbb{Z})=\bigoplus\limits_{i=0}^n{\check H_{i}(X,f;\mathbb{Z})};\\
    C_{*}(X,f;\mathbb{Z})=\bigoplus\limits_{i=0}^n{C_{i}(X,f;\mathbb{Z})};\\
    B_{*}(X,f;\mathbb{Z})=\bigoplus\limits_{i=0}^n{B_{i}(X,f;\mathbb{Z})}.
    \end{array}
 \right.$$
 \end{small}

\begin{lemma}\label{19}
Let $X$ be a Hausdorff space, $f\in C^0(X)$,
for which there is $n=n(J)\leq n(J,f)$ such that $\check{H}_{p}(X;\mathbb{Z})$ and $\check{H}_p(X,f;\mathbb{Z})$ exist,
where $p\leq n$ and $J$ is a ordered set induced by the set of all covers of $X$.
Then $f$ induced a linear transformation $f_*$ on
$\check{H}_{*}(X;\mathbb{Z})$,
$C_*(X;\mathbb{Z})$and on
$\check{H}_{*}(X,f;\mathbb{Z})$, respectively.
If $E_{f_*}$ is the set of all eigenvalue of $f_*$ and $\|E_{f_*}\|=\sup\{|a|;a\in E_{f_*}\}$,
then obtain the inequality:

(5)$\left.
   \begin{array}{l}
   \|E_{f_*}|_{\check{H}_*(X,f;\mathbb{Z})}\|\leq\|E_{f_*}|_{Z_*(X,f;\mathbb{Z})}\|\leq\|E_{f_*}|_{C_*{(X,f;\mathbb{Z})}}\|,\\
   \|E_{f_*}|_{H_*(X;\mathbb{Z})}\|\leq\|E_{f_*}|_{Z_*(X;\mathbb{Z})}\|\leq\|E_{f_*}|_{C_*{(X;\mathbb{Z})}}\|,\\
   \|E_{f_*}|_{C_*(X;\mathbb{Z})}\|\leq\|E_{f_*}|_{C_*(X,f;\mathbb{Z})}\|.
    \end{array}
 \right.$
\end{lemma}
\begin{proof}
Following, Lemma~$\ref{15},~\ref{18}$.
\end{proof}

What's more,
we can define $L_{C^0}$ category that its objects are $X^f$ and
its morphisms are continuous maps,
where $X$ is a Hausdorff spaces and $f$ is a continuous self-map on $X$.
Similarly,
we can define $\tilde{L}_{C^0}$ category that its objects are $\check{H}_*{(X,f;\mathbb{Z})}$ and
its morphisms are $F_*$,
where $X,Y$ are Hausdorff spaces,
$f\in C^0(X),g\in C^0(Y)$,
and $F_*$ is induced by the continuous map $F:X^f\to Y^g$.
Also, we can define homotopy and homeomorphism from $X^f$ to $X^g$ and research relations between the elements of $L_{C^0}$ and $\tilde{L}_{C^0}$.

\begin{definition}\label{tongpei20}
Let $X,Y$ be compact Hausdorff spaces,$f\in C^0(X),g\in C^0(Y)$.

(a)Define $X^f$ and $Y^g$ is $L_{1}$-homotopy equivalence,
if exist a pair of continuous mapping:

$\left.
   \begin{array}{ll}
    F:X^f\xrightarrow{}Y^g,&D:Y^g\xrightarrow{}X^f,\\
    F\circ D=id_{Y^g},&D\circ F=id_{X^f}.\\
    \end{array}
 \right.$

(b)Define $h,r:X^f\xrightarrow{}Y^g$ is $L_{2}$-homotopy,
if exist a continuous mapping:

$\left.
   \begin{array}{ll}
    F:X^f\times [0,1]\xrightarrow{}Y^g,&\\
    F(X^f,0)=h(X^f),&F(X^f,1)=r(X^f)\\
    \end{array}
 \right.$

so $h$ induces a homomorphism $h_*:\check{H}_*{(X,f;\mathbb{Z})}\xrightarrow{}\check{H}_*{(Y,g;\mathbb{Z})}$, and $r_*$ by $r$.
\end{definition}

Let $L$ be the class of set:
$$\{X^f;X \text{ is compact Hausdorff spaces},f\in C^0(X)\},$$
for each $X^f,Y^g\in L$, let $mor_s(X^f,Y^g)=L_{1}(X^f,Y^g)$,
by the $L_{1}$-homotopy and composition of function $\circ$, we get a category$(L,mor_s,\circ)$.

Let $\tilde{L}$ be the class of set:
$$\{\check{H}_*{(X,f;\mathbb{Z})};X^f\in L \},$$
for $\check{H}_*{(X,f;\mathbb{Z})},\check{H}_*{(Y,g;\mathbb{Z})}\in\tilde{L}$,
let $mor_H(\check{H}_*{(X,f;\mathbb{Z})},\check{H}_*{(Y,g;\mathbb{Z})})$ be the all group homomorphisms from $\check{H}_*{(X,f;\mathbb{Z})}$ to $\check{H}_*{(Y,g;\mathbb{Z})}$,
by the $L_{1}$-homotopy induced $*$ mapping and composition of function $\circ$, we get a category$(\tilde{L},mor_H,\circ)$.
Easy to see the functor from $(L,mor_s,\circ)$ to $(\tilde{L},mor_H,\circ)$.

\begin{theorem}\label{tonglun21}
Let $f\in C^0(X),g\in C^0(Y)$, $X,Y$ becompact Hausdorff spaces,

(a)if $X^f$ and $Y^g$ is $L_{1}$-homotopy equivalence,

then ${C}_{p}(X,f;\mathbb{Z})={C}_{p}(X,g;\mathbb{Z})$ and $\check{H}_{p}(X,f;\mathbb{Z})=\check{H}_{p}(X,g;\mathbb{Z})$.

(b)if $h,r:X^f\xrightarrow{}Y^g$ is $L_{2}$-homotopy,

then $h_*=r_*$.
\end{theorem}
\begin{proof}
By diagram chasing.
\end{proof}

\begin{example}
Let $f\in C^0(X),g\in C^0(Y)$, $X,Y$ be compact Hausdorff spaces,
if there is a homeomorphism $F$ from $X$ to $Y$ such that $Ff=gF$,
then $$\check{H}_{p}(X,f;\mathbb{Z})=\check{H}_{p}(X,g;\mathbb{Z}).$$
\end{example}

\begin{example}
Let $f\in C^0(X),g\in C^0(Y)$, $X,Y$ be compact Hausdorff spaces,
if there is a homeomorphism $F$ from $X$ to $Y$,
then $$\check{H}_{p}(X,f;\mathbb{Z})=\check{H}_{p}(X,g;\mathbb{Z}).$$
\end{example}

\begin{example}
Let $f\in C^0(X),g\in C^0(Y)$, $X,Y$ be compact Hausdorff spaces,
if there is a continuous mapping $F:X\times [0,1]\xrightarrow{}Y$ such that:

$\left\{
   \begin{array}{l}
    F(X,0)=h(X),\\
    F(X,1)=r(X),
    \end{array}
 \right.$

that is,  $h$ and $r$ are homotopy.
Then $h_*=r_*$,
where
$$h_*:\check{H}_*{(X,f;\mathbb{Z})}\xrightarrow{}\check{H}_*{(Y,g;\mathbb{Z})},$$
$$r_*:\check{H}_*{(X,f;\mathbb{Z})}\xrightarrow{}\check{H}_*{(Y,g;\mathbb{Z})}.$$
\end{example}

\section{Topological Entropy Conjecture}

In this part,
$X$ is a compact Hausdorff space and $J$ is the set of all covers
for which there is $n(J)$ such that $\check{H}_*(X;\mathbb{Z})$ exists.
For $\alpha,\beta$ are collections of open sets of $X$,

define:

(6)$\left\{
   \begin{array}{ll}
    \alpha\vee\beta=\{A\cap B;A\in\alpha,B\in\beta\};\\
    f^{-1}(\alpha)=\{f^{-1}(A);A\in\alpha\},\\
    f^{-1}(\alpha\vee\beta)=f^{-1}(\alpha)\vee f^{-1}(\beta);\\
    \bigvee\limits_{i=0}^{n-1} f^{-i}(\alpha)=\alpha\vee f^{-1}(\alpha)\vee\cdots\vee f^{-(n-1)}(\alpha).
    \end{array}
 \right.$

If $\alpha,\beta$ are open covers of $X$,
let $N(\alpha)$ is the infimum numbers of subcover of $\alpha$.
For the compact of $X$,
we know that $N(\alpha)$ is a positive integer.
so define

$$H(\alpha)=\log N(\alpha)\geq0.$$

Following \cite{PLPeterWalters}P81, there is

$$\alpha<\beta\Longrightarrow H(\alpha)\leq H(\beta).$$

\begin{definition}[\cite{PLPeterWalters}P89]\label{22}
For a fix open cover $\alpha$ of $X$,define :

$$ent(f,\alpha)=\lim\limits_{n\rightarrow\propto}\frac{1}{n}H(\bigvee\limits_{i=0}^{n-1} f^{-i}(\alpha)),$$

and define the topological entropy of $f$:

$$ent(f)=\sup\limits_{\alpha}\{ent(f,\alpha)\},$$

where $\sup\limits_\alpha$ is through the all open covers of $X$.
\end{definition}

Let $\alpha$ be a open cover of $X$ and $L_{f}(\alpha)=\{L_{f}(U)|U\in\alpha\}$,
then $L_{f}(\alpha)$ can be induced a open fiber cover $\dot L_{f}(\alpha)$ of $X^f$.

\begin{definition}\label{24}
For a fix open fiber cover $\dot L_f(\alpha)$ of $X^f$,
define :

$\left\{
   \begin{array}{ll}
\frac{f^{-1}(\dot L_f(\alpha))}{\dot L_f(\alpha)}=\max\limits_{U\in \alpha}\|\{f^{-1}\dot L_f(U)\cap \dot L_f(U)\}\|;\\
\frac{f(\dot L_f(\alpha))}{\dot L_f(\alpha)}=\max\limits_{U\in \alpha}\|\{f\dot L_f(U)\cap \dot L_f(U)\}\|;\\
L_d=\max\{\frac{f^{-1}(\dot L_f(\alpha))}{\dot L_f(\alpha)},\frac{f(\dot L_f(\alpha))}{\dot L_f(\alpha)}\};\\
ent(f_L,\dot L_f(\alpha))=ent(f,\alpha)+\log L_d.
    \end{array}
 \right.$

and define the topological fiber entropy of $f$:

$$ent(f_L)=\sup\limits_{\dot L_f(\alpha)}\{ent(f_L,\dot L_f(\alpha))\},$$

where $\sup\limits_{\dot L_f(\alpha)}$ is through the all open covers of $X^f$.
\end{definition}

\begin{lemma}[\cite{PLPeterWalters}P102]\label{25}
If $f$ is the shift operator on a $k$-symbolic space, then $ent(f)=\log k$.
\end{lemma}

\begin{corollary}\label{26}
If $f$ is the shift operator on a $k$-symbolic space, then
$$ent(f_L)=ent(f)+\log k=2\log k.$$
\end{corollary}

\begin{example}
Let $\{1,2,\cdots,k\}=X$,

if define
$f:\left\{
   \begin{array}{ll}
   \{1\}\rightarrow\{1,2,\cdots,k\},\\
   \{2\}\rightarrow\{1,2,\cdots,k\},\\
   \quad\vdots\quad \vdots\qquad\vdots\\
   \{k\}\rightarrow\{1,2,\cdots,k\}
    \end{array}
 \right.$,

then $ent(f)=0,ent(f_L)=\log k$.
\end{example}

\begin{example}
Let $\{1,2,\cdots,k\}=X$,
if define
$f:\left.
   \begin{array}{ll}
   \{1,2,\cdots,k\}\to\{1\}
    \end{array}
 \right.$

then $ent(f)=0,ent(f_L)=\log k$.
\end{example}

\begin{example}
Let $[0,1]=X$ ,if define $f(x)=kx,   0<k<1$.

then $ent(f)=0,ent(f_L)=-\log k$.
\end{example}

\begin{lemma}\label{27}
Let $2< m\in\mathbb{Z}$,
then there are $p,q\in \mathbb{Z}$ such that
$p\neq q,m=p+q$, where $1<p,1<q$.
\end{lemma}

Let $f\in C^0(X)$,
$f_*$ be the linear transformation on $\check{H}_*(X,f;\mathbb{Z})$ induced by $f$.
Define \v{C}ech eigenvalue chains is the chains belong to the eigenvalues of $f_*$,
and any \v{C}ech eigenchains can be induced a open cover of $X^f$.

\begin{lemma}\label{28}
Let $X$ be a compact Hausdorff space
for which there is $n=n(J)\leq n(J,f)$ such that $\check{H}_{p}(X;\mathbb{Z})$ and $\check{H}_p(X,f;\mathbb{Z})$ exist,
where $p\leq n$ and $J$ is a ordered set induced by the set of all covers of $X$.
Let $\alpha\in J$ be a open cover,
if $L_{f}(\alpha)$ is a $\check{C}ech$ eigenchains belong to eigenvalue $m$,
then $L_{f}(\alpha)$ has a factor which conjunct to a shift operator on $m$-symbolic space or its $L_d=m$,
where $0\leq m\in \mathbb{Z}$.
\end{lemma}
\begin{proof}
Let $L_f=\sum\limits_{i=0}^{k}{a_i}\check\sigma_i$
be a eigenchains belong to the eigenvalue $m$,
where $\check\sigma_i\in \check{H}_*(X,f;\mathbb{Z})$, $m,{a_i}\in\mathbb{Z}$.

From Lemma~$\ref{18}$
there is $f$-\v{C}ech homology germ $H_{p}(J,f;\mathbb{Z})$ such that
$$H_{p}(J,f;\mathbb{Z})\sim\check H_{p}(X,f;\mathbb{Z}),0\leq p\leq n(J).$$

Hence, there is $\Phi\in J$ such that $L_f\in H_{*}(\Phi,f;G)$ and
$$f_*(L_f)=m(L_f).$$

That can be extend to a equation on $C_{*}(\Phi,f;G)$,
and get the equation

$$f_\sharp(\check\sigma_i)=m(\check\sigma_i),i\in\{0,\cdots,k\},$$
where $\check\sigma_i\in C_*(\Phi,f;G)$ and $m,{a_i}\in \mathbb{Z}$.

Just thinking $f_\sharp$ on $C_{*}(\Phi,f;G)$,
let $U_0,\cdots,U_j$ is subset of $X$ and
$$\check\sigma_i=L_f(U_0)\cap\cdots\cap L_f(U_j)$$ then

$$L_{f}(U_{\eta})=\{\cdots,{f}^{-n}(U_{\eta}),\cdots,{f}^{-1}(U_{\eta}),{f}^{0}(U_{\eta}),{f}^{1}(U_{\eta})\cdots,{f}^{n}(U_{\eta}),\cdots\},$$
where $\eta\in\{0,\cdots,j\}.$

So
$$\left.
   \begin{array}{ll}
f_\sharp(\check\sigma_i)=f_\sharp(L_f(U_0)\bigcap\cdots\bigcap L_f(U_j))\\
=L_{f}(f(U_0))\bigcap\cdots\bigcap L_{f}(f(U_j))\\
=m(L_f(U_0)\bigcap\cdots\bigcap L_f(U_j)).
    \end{array}
 \right.$$

That is,

\begin{small}
$$\left.
   \begin{array}{ll}
   m(\bigcap\limits_{\eta=0}^{j}
\{\cdots,{f}^{-n}(U_\eta),\cdots,{f}^{-1}(U_\eta),\underline{{f}^{0}(U_\eta)},\cdots,{f}^{n}(U_\eta),\cdots\})\\
=\bigcap\limits_{\eta=0}^{j}\{\cdots,{f}^{-n}(f(U_\eta)),\cdots,{f}^{-1}(f(U_\eta)),\underline{{f}^{0}(f(U_\eta))},\cdots,{f}^{n}(f(U_\eta)),\cdots\}\\
=\bigcap\limits_{\eta=0}^{j}\{\cdots,{f}^{-(n-1)}(f(U_\eta)),\cdots,{f}^{-1}(f(U_\eta)),\underline{f(U_\eta)}, f^2(U_\eta),\cdots,{f}^{n+1}(U_\eta),\cdots\}.
    \end{array}
 \right.$$
 \end{small}

Therefore $$m(\bigcap\limits_{\eta=0}^{j}L_f(U_\eta))=(\bigcap\limits_{\eta=0}^{j}L_f(f(U_\eta))).$$

Loss no generally let $j=0$,
then
$$L_f(f(U_0))=m(L_f(U_0)).$$

We only prove this for torsion free $L_f(U_0)$,
for torsion $L_f(U_0)$ the conclusion is trivial.
Now let $L_f(U_0)$ be a torsion free element.

(i) $m=0,1$ the conclusion is trivial.

(ii) $m=2$, exist $U\subseteq f^{-1}(f(U_0))$,$U\nsubseteq U_0 $ and $U_0\nsubseteq U$,
$U_0,U$ are non-empty open subsets.

else $f^{-1}(f(U_0))=U_0$ and there is

$$L_f(f(U_0))=(L_f(U_0))=2(L_f(U_0)),$$

this is a contradiction for the properties that $\mathbb{Z}$ is a free group.

Because of $U\nsubseteq U_0 $ and $U_0\nsubseteq U$,$U_0,U$ is open subsets,
with the property of Hausdorff space,
there are points $x\in U_0,x\notin U,y\in U,y\notin U_0$ and open neighborhoods $O(x)$ and $O(y)$,
respectively,
such that $x\in O(x)\in U_0,O(x)\notin U$ and $y\in O(y)\in U,O(y)\notin U_0$.
Then $O(x),O(y)\subseteq f^{-1}(f(U_0))$,
and $O(x)\bigcap O(y)=\emptyset$.

so $L_d=2$, and for $m=2$ the conclusion is true.

(iii) $m\geq3$, by induction, for $m=n-1$ the conclusion is right, look $m=n$,

using $lemma~\ref{27}$, we get $m=p+q,p\neq q$ and
$$L_f(f(U_0))=p(L_f(U_0))+q(L_f(U_0)).$$

Therefore there is $f|_{U_0}=h+g$ such that
$$L_h(f(U_0))=p(L_f(U_0)),L_g(f(U_0))=q(L_f(U_0)).$$

If $L_f\neq L_{h\oplus g}$,
then using $(ii)$ with the same computing we get $L_d=m$.

If $L_f= L_{h\oplus g}$,
then there is
$$L_f(f(U_0))=L_h(f(U_0))\oplus L_g(f(U_0)),$$

else there is
$$h^{-1}(f(U_0))\bigcap g^{-1}(f(U_0))=W\neq\emptyset,p(L_f(W))=q(L_f(W))$$,
but $p\neq q$,
contradict with the property that $\mathbb{Z}$ is a free group.

For $m=p+q$ we get that $p,q\leq n-1$, by the induction there is

$$\left\{
   \begin{array}{ll}
    h^{-1}(f(U_0))\supseteq U_{0i},U_{0j},& U_{0i}\bigcap U_{0j}=\emptyset,1\leq i,j\leq p,\\
    g^{-1}(f(U_0))\supseteq U_{1k},U_{1l} & U_{1k}\bigcap U_{1l}=\emptyset,1\leq k,l\leq q,\\
    \end{array}
 \right.$$

where $U_{0i},U_{0j},U_{1k},U_{1l}$ are non-empty open subset.

With the decomposition

$$L_f(f(U_0))=L_h(f(U_0))\oplus L_g(f(U_0)),$$

we get that
$f^{-1}(f(U_0))\supseteq U_i,U_j$,
$U_i\bigcap U_j=\emptyset,1\leq i,j\leq m$ and they all are non-empty open subset of $X$.

Therefore we can get a $m$-symbolic space $S_m$,
that is $L_{f}(U_0)$ has a factor conjunct a shift operator on $S_m$.

So for $m=n$ the conclusion is right,
by the induction the conclusion is right for any eigenvalue $m,0\leq m\in \mathbb{Z}$.

\end{proof}

\begin{theorem}\label{29}
Let $X$ be a compact Hausdorff space
for which there is $n=n(J)$ such that $\check{H}_{p}(X;\mathbb{Z})$ exist,
where $p\leq n$ and $J$ is a ordered set induced by the set of all covers of $X$.
For $f\in C^0(X)$, there is
$$\log\|E_{f_*{\check{H}_*(X;\mathbb{Z})}}\|\leq ent(f_L),$$

Moreover, if $\check{H}_p(X,f;\mathbb{Z})$ exist for $p\leq n(J,f)$,
then there is
$$\log\|E_{f_*{\check{H}_*(X,f;\mathbb{Z})}}\|\leq ent(f_L).$$
\end{theorem}

\begin{proof}
Following Lemma$~\ref{15},~\ref{19},~\ref{28}$,
we get:

$$ent(f_L)\geq ent(f_L,\dot L_(\alpha))\geq\log \| E_{f_*}|_{C_*(X,f;\mathbb{Z})}\|
\geq\log \| E_{f_*}|_{\check{H}_*(X;\mathbb{Z})}\|,$$
and

$$ent(f_L)\geq ent(f_L,\dot L_(\alpha))\geq\log \| E_{f_*}|_{C_*(X,f;\mathbb{Z})}\|
\geq\log \| E_{f_*}|_{\check{H}_*(X,f;\mathbb{Z})}\|.$$
\end{proof}

For simple computing we get the following
\begin{proposition}\label{gongebubiandeqianweishang}
$ent(I_L)$ is invariant under topological conjugate.
\end{proposition}

\begin{proposition}\label{30}
$ent(f_L)\geq ent(f)$,inequality can be strict.
\end{proposition}

\begin{proposition}\label{31}
$ent(I_L)=ent(I)=0$, where $I$ is identical mapping.
\end{proposition}

\begin{corollary}\label{32}
Let $X$ be a compact Poincar\`{e} space
for which there is $n=n(J)$ such that $\check{H}_p(X;\mathbb{Z})$ exists,
where $p\leq n$ and $J$ is a ordered set induced by the set of all covers of $X$.
Topological Entropy Conjecture for topological fiber entropy is valid for \v{C}ech cohomology.
Moreover if $\check{H}_p(X,f;\mathbb{Z})$ exists for $p\leq n(J,f)$,
then Topological Entropy Conjecture for topological fiber entropy is valid for $f$-\v{C}ech homology.
\end{corollary}

\begin{corollary}\label{34}
In triangulable compact $n$-dimension manifold
Topological Entropy Conjecture is valid for topological fiber entropy and homology group $H_{*}(M;\mathbb{Z})=\bigoplus\limits_{i=0}^n{H_{i}(M;\mathbb{Z})}$,
where $H_{i}(M;\mathbb{Z})$ is the $i$th integer coefficients homology group of $M$.
\end{corollary}

\begin{remark}
If we replace $\mathbb{Z}$ with any free abelian group $G$ of finite spanning,
then the conclusion is also valid  .
\end{remark}

\begin{remark}
Because the counterexample of Katok \cite{ABKatok1986} is on $2$-dimension sphere $S^{2}$ and $f\in C^{0}(S^2)$,
with corollary~$\ref{34}$ we get that the conjecture is valid again with our definition,
that is,
$\log\rho\leq ent(f_{L})$.
\end{remark}

\begin{remark}
To establish the inequality $\log\rho\leq ent(f_{L})$ is my interesting,
but others maybe more interest what does topological fiber entropy $ent(f_{L})$ measure?
From the definition
$$ent(f_L)=\sup\limits_{\dot L_f(\alpha)}\{ent(f,\alpha)+\log L_d\},$$
we get that topological fiber entropy $ent(f_{L})$ is make $\sup\limits_{\dot L_f(\alpha)}$ on the sum $$ent(f,\alpha)+\log L_d.$$
The first part $ent(f,\alpha)$ is the usually one.
But the second part $\log L_d$ is likely some fiber ratio of the dynamics $(X,f)$,
it is likely the `reference system' or `initial value' of the first part $ent(f,\alpha)$.
\end{remark}

\section{Acknowledgments:}
The author would like to thank the referee for his/her careful reading of the paper and helpful comments and suggestions.
I appreciate all the experiences that make me grow up.
Thank you for all the people in our life.
Years past, my English is still accompany my weakness,
but, the latest version is more like a paper than the old one.
Any way, the old version have more inspired.
Recently, I am in the vein for reversion this paper,
and it is too long time to keep some thing well-founded,
such as, I can't find where is cited the reference \cite{CP,AM} within this paper,
and just keep them here as the original version.



\begin{thebibliography}{99}\small
\bibitem{FrankWAndersonKentRFuller1992}Frank W. Anderson and Kent R. Fuller, \textit{Rings and Categories of Modules}, Springer-Verlag,2nd, New York, USA, 1992.

\bibitem{ABKatok1986}A. B. Katok, \textit{A conjecture about Entropy}, Amer. Math. Soc. Transl(2)Vol. 133,1986.

\bibitem{JamesRMunkres2006}James R. Munkres, \textit{Elements of Algebraic Topology}, Science Press, Peking, China, 2006.

\bibitem{PLPeterWalters}Peter Walters, \textit{An Introduction to Ergodic Theory}, Springer-Verlag New York,Inc.,1982.
\bibitem{MShub1974}M. Shub, \textit{Dynamical System,fitrations and entropy}, Bull. Amer. Math. Soc. \small{\textbf{80}}(1974), 27-41.

\bibitem{BLvanderWaerden1991}B. L. van der Waerden, \textit{Algebra}, Springer-Verlag New York,Inc.,1991

\bibitem{CP}Charles Pugh, \textit{On the entrony conjecture: a report on conversations among R. bowen, M. Hirsch, A. Manning, C. Pugh, B. Sanderson, M. Shub, and R. Williams}, in$^{[6]}$,pp,257-261.

\bibitem{AM}Anthony Manning{editor}, \textit{Dymnamical systems-Warwick1974 (Proc. Sympos., Coventry, 1973/1974, Presented to E. C. Zeeman on His Fiftieth Birthday)}, Lecture Notes in Math, vol.468, Springger-Verlag, 1975.

\end{thebibliography}
\end{document}